\documentclass[11pt]{article}

\usepackage{graphicx,graphics}
\usepackage{amsmath,amssymb,amsthm,graphics,amsfonts,graphicx}
\usepackage{appendix}
\usepackage{lineno}
\newtheorem{theorem}{Theorem}[section]

\begin{document}
\begin{linenumbers}

\begin{center}  {\Large \bf Analysis of a malaria model with mosquito host choice and bed--net control}
\end{center}
\smallskip
\begin{center}
{\small \textsc{Bruno Buonomo}}
\end{center}

\begin{center} {\small \sl Department of Mathematics and Applications,
University of Naples Federico II\\  via Cintia, I-80126 Naples,
Italy\\ buonomo@unina.it}\\
\end{center}

%
{\centerline{\bf Abstract}
\begin{quote}
\small  A malaria model is formulated which includes the enhanced attractiveness of infectious humans to mosquitoes,  as result of host manipulation by malaria parasite, and the human behavior, represented by insecticide-treated bed nets usage. The occurrence of a backward bifurcation at $R_0=1$ is shown to be possible, which implies that multiple endemic equilibria co-exist with a stable disease-free equilibrium when the basic reproduction number is less than unity. This phenomenon is found to be caused by disease--induced human mortality. The global asymptotic stability of the endemic equilibrium for $R_0>1$ is proved, by using the geometric approach to global stability. Therefore, the disease becomes endemic for $R_0>1$ regardless of the number of initial cases in both the human and vector populations. Finally, the impact of vector's host preferences and bed--net usage behavior on system dynamics is investigated.
\end{quote}}
%

\noindent { {\bf Subject class}: 92D30, 34C23, 34D23}




\section{Introduction} 
\label{intro}

Malaria is a life--threatening disease caused by parasites transmitted to susceptible humans through the bites of infected female mosquitoes of the genus \textit{Anopheles}. In spite of recent successes in the struggle against malaria that have lead to a substantial reduction of reported malaria cases and deaths, latest estimates indicate that malaria is still a global emergency with 219 million cases in 2010 and a death toll ranging from 660 000 individuals \cite{who} to 1,24 million \cite{muetal}.

Mathematical modeling of malaria transmission has a long history, intimately linked to the evolution and history of malaria over more than 100 years. As part of the necessary multi--disciplinary research approach, mathematical models have been used to provide a framework for understanding malaria transmission dynamics and the best strategies to control the disease.
Starting from the basic Ross--MacDonald models \cite{mac,ro}, a very large literature on the subject is nowadays available. A comprehensive survey on malaria mathematical modeling may be got from classical sources, e. g. \cite{arma,ko,ne}, and more recent contributions \cite{masasi,smetal}.

In this paper, we focus on two specific aspects that have recently received much attention from malaria modelers: \emph{(a)} the enhanced attractiveness of infectious humans to mosquitoes; \emph{(b)} the non-pharmaceutical interventions (NPI) usage.

As for point \emph{(a)}, it concerns the investigation of behavioral manipulation by malaria parasite to increase the host's attractiveness to mosquitoes \cite{pole}. In this direction, experimental evidence has revealed the enhanced attractiveness to mosquitoes of hosts harboring the parasite's gametocytes  (the stage infective to mosquitoes) \cite{laetal}. 
Several mathematical models have been proposed to get an insight on understanding and prediction of disease evolution when  \emph{vector bias} to infected hosts is taken into account \cite{bbcruz2,chbr,hoetal,ki,vdl}. In particular, the model introduced in \cite{chbr} is obtained by extending the classical Ross model \cite{ro} to include the enhanced attractiveness of infectious humans to mosquitoes. Later, this model has been further extended to include both immigration and disease--induced death of humans \cite{bbcruz2}.

As for point \emph{(b)}, non-pharmaceutical interventions aim to limit virus spread by reducing contact between infectious and susceptible individuals \cite{linetal}. The insecticide--treated bed--nets (ITN) are among the NPI specifically targeted for malaria transmission \cite{le}. The effectiveness of ITN is largely influenced by behavioral factors. In fact, peoples may decide to not adopt ITN, in spite of its usefulness \cite{le}, because of personal reasons. Especially during the dry season, hot weather, tendency to sleep outside the house and lack of mosquito nuisance are among the reasons for not using the ITN \cite{fretal}. As a consequence, the role of human behavior (and misbehaviors) ought be included in the modeling of ITN--usage. In this case, modeling in the framework of \emph{Behavioral Epidemiology}, where the key aspect is the impact of human behavior on epidemics, is appropriate \cite{thebook}.

In this paper, we formulate and analyze a malaria model that includes both vector--bias preference for infectious host and bed--net usage from the population. The baseline model is the vector--bias malaria model proposed in \cite{bbcruz2}, which we extend by adopting the bed--net usage modeling proposed in \cite{agetal}, i.e. it is assumed that the contact rate and the mosquito mortality are functions of the bed--net usage.

We perform a bifurcation analysis to detect the occurrence of a backward bifurcation and, consequently, the presence of multiple endemic equilibria co-existing with a stable disease-free equilibrium when the basic reproduction number is less than unity.  This result is obtained by using the centre manifold theory \cite{ccso,duhucc,vawa02}. Under the point of view of disease control, the occurrence of backward bifurcation has very important implications because the classical threshold condition, $R_0<1$, is no longer sufficient to obtain the elimination of the disease from the population.

We also perform a global stability analysis of the endemic equilibrium for $R_0>1$, in the case that the total vector population is at equilibrium.  We use the geometric method to stability \cite{limu}.  This result ensures that the disease will become endemic for $R_0>1$ whatever the initial cases of infection in both the populations are.

Finally, we assess both the individual and simultaneous impact of  bed--net usage behavior and vector--bias preferences on system dynamics and, in particular, how they can influence the basic reproduction number and the occurrence of the backward bifurcation.

The rest of the paper is organized as follows: in Section \ref{sec:model} we introduce the model and give some basic properties, including the local stability of the disease--free equilibrium and the existence of endemic states. In Section \ref{bifurc} we perform the bifurcation analysis. Section \ref{global} is devoted to the global stability analysis of the endemic state. In Section \ref{impact} the impact of bed--net usage behavior and vector--bias preferences on system dynamics is discussed. Concluding remarks are given in Section \ref{concl}.

\section{Model and basic properties} 
\label{sec:model}

We consider the following system of nonlinear ordinary differential equations:
\begin{equation}\label{model}
\begin{array}{ll}
\dot S_{h}=&\Lambda_h-\lambda_h(b,\pi)S_h-\mu S_{h}+\delta I_{h}\\
\dot I_{h}=&\lambda_h(b,\pi)S_h-(\alpha+\mu+\delta) I_{h}\\
 \dot S_{v} =&\Lambda_v-\lambda_v(b,\pi)S_{v}-\eta(b) S_{v}\\
 \dot I_{v}=& \lambda_v(b,\pi)S_{v}-\eta(b) I_{v},
\end{array}
\end{equation}
where the upper dot denotes the time derivative and the state variables are given by susceptible humans, $S_h$, infectious humans, $I_h$, susceptible vectors, $S_v$ and infectious vectors $I_v$. The parameters are all strictly positive constants and their meaning is described in Table \ref{Tab:param}. The terms $\lambda_h(b,\pi)$ and $\lambda_v(b,\pi)$ are the forces of infection.

\begin{table}
\caption{Description and baseline values of parameters in system (\ref{model}).}
\label{Tab:param}       
\begin{tabular}{lll}
\hline\noalign{\smallskip}
Parameter & Description & Baseline value \\
\noalign{\smallskip}\hline\noalign{\smallskip}
$\Lambda_h$ & Immigration rate in humans &  $10^3/(70\times365)$  \\ 
$\Lambda_v$ & Immigration rate in mosquitoes & $10^4/21$\\
$\pi$ & vector--bias parameter & varies\\
$b$ & Proportion of ITN usage & varies \\
$\mu$ & Natural mortality rate in humans & $1/(70\times365)$\\
$\nu_{nat}$ & Natural mortality rate in mosquitoes & $1/21$\\
$\nu_{bn}$ & maximum NTI-induced death rate in mosquitoes & $1/21$\\
$\alpha$ & Disease--induced death rate in humans & $10^{-3}$\\
$p_1$ & Prob. of disease transm. from mosquito to human& 1\\
$p_2$ & Prob. of disease transm. from human to mosquito & 1\\
$\beta_{\max}$ & Maximum transmission rate& 0.1\\
$\beta_{\min}$ & Minimum transmission rate& 0\\
$\delta$ & Recovery rate of infectious humans to be susceptible& 1/4\\ 
\noalign{\smallskip}\hline
\end{tabular}
\end{table}

Following \cite{chbr} we assume that mosquitoes that bite humans will do it at probability $p$ if the human is infectious,
and probability $q$, with $p>q$, if the human is susceptible. Hence, when a mosquito bites an human, the probability that this
human is infected is given by the ratio between the total bitten infectious humans and the total bitten humans, $pI_h/(pI_h+qS_h)$, whereas the probability that this human is susceptible is given by the ratio between the total bitten susceptible humans and the total bitten humans, $qS_h/(pI_h+qS_h)$.

If $\pi$ denote the ratio $p/q$, then $\pi\geq  1$ ($\pi=1$ means that the enhanced attractiveness of infectious humans to mosquitoes is neglected) and  it follows:
\begin{equation}\label{lambda}
\lambda_h(b,\pi)=\frac{p_1\beta(b)\,I_v}{\pi I_h + S_h},\;\;\;\;\;\;
\lambda_v(b,\pi)=\frac{\pi p_2\beta(b)\,I_h}{\pi I_h + S_h}.
\end{equation}

The parameter $b$ represents the bed net usage. It ranges between $b=0$ (no bed net users) and $b=1$ (all the individuals of host population are users). Using bed nets reduces the probability for humans to be bitten. Moreover, the nets are treated with insecticide. Therefore  the role of $b$ in the model is to reduce the contact rate $\beta$ and to increase the mosquito death rate $\eta$. Therefore, it is assumed that \cite{agetal}
\begin{equation}\label{etab}
\eta(b)=\eta_{nat}+\eta_{bn}b,\;\;\;\;\beta(b)=\beta_{\max}-b\left(\beta_{\max}-\beta_{\min}\right),\;\;\;\;\;0\leq b\leq1.
\end{equation}

Denote by $N_h$ and $N_v$ the total human and vector population, respectively (i.e. $N_h=S_h+I_h$ and $N_v=S_v+I_v$).  Note that
\begin{equation}\label{TotPop}
\dot N_h = \Lambda_h - \mu N_h - \alpha I,\;\;\;\;\;\dot N_v = \Lambda_v - \eta(b) N_v.
\end{equation}
From these equalities, by using positiveness of solutions and a comparison theorem \cite{laetal},  it is not difficult to show (see \cite{agetal}) that model \eqref{model} can be studied in the positively invariant and attractive set
$$
\Omega=\left\{(S_h,I_h,S_v,I_v)\in{\bf R}^4: 0\leq N_h(t) \leq\frac{\Lambda_h}{\mu},\;\;0\leq N_v(t)\leq\frac{\Lambda_v}{\eta(b)}\right\}.
$$
System \eqref{model} admits the \emph{disease-free equilibrium}
\begin{equation}\label{E0}
E_0:=\left(S_{h0},\, 0,\, S_{v0},\, 0\right)=\left(\frac{\Lambda_h}{\mu},\, 0,\, \frac{\Lambda_v}{\eta(b)},\, 0\right).
\end{equation}
The Jacobian matrix corresponding to \eqref{model} is,
\begin{equation}\label{Jacobian}
J=
\left(\begin{array}{cccc}
-\frac{\partial \lambda_h}{\partial S_h}S_h-\lambda_h-\mu & -\frac{\partial \lambda_h}{\partial I_h}S_h + \delta & 0 &  -\frac{\partial \lambda_h}{\partial I_v}S_h\\
\frac{\partial \lambda_h}{\partial S_h}S_h+\lambda_h& \frac{\partial \lambda_h}{\partial I_h}S_h-\alpha_0& 0 & \frac{\partial \lambda_h}{\partial I_v}S_h\\
-\frac{\partial \lambda_v}{\partial S_h}S_v& -\frac{\partial \lambda_v}{\partial I_h}S_v& -\lambda_v-\eta(b) & 0 \\
\frac{\partial \lambda_v}{\partial S_h}S_v & \frac{\partial \lambda_v}{\partial I_h}S_v & \lambda_v& -\eta(b)
\end{array}
\right),
 \end{equation}
 where $\alpha_0=\alpha+\mu+\delta$, and, in view of \eqref{lambda},
 $$
 \frac{\partial \lambda_h}{\partial S_h}=-\frac{p_1\beta(b)\,I_v}{\left(\pi I_h+S_h\right)^2};\;\;\;\;\;
  \frac{\partial \lambda_h}{\partial I_h}=-\frac{p_1\pi\beta(b)\,I_v}{\left(\pi I_h+S_h\right)^2};\;\;\;\;\;
  \frac{\partial \lambda_h}{\partial I_v}=\frac{p_1\beta(b)}{\pi I_h+S_h},
 $$
 and
  $$
 \frac{\partial \lambda_v}{\partial S_h}=-\frac{\pi p_2\beta(b)\,I_h}{\left(\pi I_h+S_h\right)^2};\;\;\;\;\;
  \frac{\partial \lambda_v}{\partial I_h}=\frac{\pi p_2\beta(b)\,S_h}{\left(\pi I_h+S_h\right)^2}.
 $$
Introduce now the basic reproduction number
\begin{equation}\label{R0}
R_0 = \frac{\pi p_1p_2\mu\Lambda_v\beta^2(b)}{\Lambda_h\eta^2(b)\left(\alpha+\mu+\delta\right)}.
\end{equation}
From now on we will omit, when it is not necessary, to explicitly indicate the $b$--dependence of $\beta$ and $\eta$. We have the following result.
\begin{theorem} \label{stabE0}
The disease--free equilibrium $E_0$, given by \eqref{E0}, is locally asymptotically stable if $R_0<1$ and unstable if $R_0>1$.
\end{theorem}
\noindent{\bf Proof}. Evaluated at $E_0$, the Jacobian matrix (\ref{Jacobian}) gives
\begin{equation}
J(E_0)=
\left(\begin{array}{cccc}
-\mu & \delta & 0 & -p_1\beta \\
0 & -\alpha_0 & 0 &  p_1\beta \\
0 & -\varphi & -\eta & 0 \\
0 &  \varphi & 0 & -\eta
\end{array}
\right),
 \end{equation}
 where
 $$
 \varphi=\frac{\pi p_2\beta\mu\Lambda_v}{\eta\Lambda_h}.
 $$
 The eigenvalues are given by $\lambda_1=-\mu$, $\lambda_2=-\eta$ and the other two are eigenvalues of the submatrix
$$\overline{J}(E_0)=\begin{pmatrix}
  -\alpha_0  & p_1\beta \\ \ \\
  \varphi & -\eta\
 \end{pmatrix}.$$
The trace of $\overline{J}$ is negative, and the determinant is
$$
\det \overline{J}(E_0) = \alpha_0 \eta\left(1-\frac{p_1\beta\varphi}{\eta\alpha_0}\right).
$$
In view of \eqref{R0} it follows $\det \overline{J}(E_0)=\alpha_0 \eta\left(1-R_0\right)$, so that $E_0$ is stable if $R_0<1$, and unstable if $R_0>1$. \#

Now let us introduce the quantities
\begin{equation}\label{A0}
A_0=\eta\Lambda_h\pi^2(\eta+p_2\beta),
\end{equation}
\begin{equation}\label{B0}
B_0=\pi\left[\eta\alpha_0\Lambda_h\left(2\eta+p_2\beta\right)-p_1p_2\beta^2\Lambda_v\left(\alpha+\mu\right)\right],
\end{equation}
\begin{equation}\label{C0}
C_0=\eta^2\alpha_0^2\Lambda_h\left(1-R_0\right).
\end{equation}
Note that $C_0>0$ is equivalent to $R_0<1$, and the $B_0<0$ is equivalent to $R_0>R_a$, where
$$
R_a=\frac{\pi\mu(2\eta+p_2\beta)}{\eta(\alpha+\mu)}.
$$
The following theorem concerns the existence of endemic equilibria:
\begin{theorem}\label{exendemic}
Model (\ref{model}) has
\begin{itemize}
\item[\emph{(i)}] a unique endemic equilibrium if $C_0<0$ (i.e. if $R_0>1$);
\item[\emph{(ii)}] a unique endemic equilibrium if
\begin{equation}
B_0<0,\;\;\;\;\; {\rm and}\;\;\; C_0=0,\;\;\;{\rm
or}\;\;\;B_0^2-4A_0C_0=0;
\end{equation}
\item[\emph{(iii)}] two endemic equilibria if
\begin{equation}
C_0>0, \;\;\;\;\; B_0<0 \;\;\;\;\; {\rm and}\;\;\;B_0^2-4A_0C_0>0;
\end{equation}
\item[\emph{(iv)}] no endemic equilibria otherwise.
\end{itemize}
\end{theorem}

\noindent {\bf Proof.} Denote by $E^*=\left(S_h^*,I_h^*,S^*_v,I^*_v\right)$ a generic endemic equilibrium of model (\ref{model}). In view of (\ref{model}), the components must be solutions of the following equations:
$$
S_h^*=\frac{\alpha_0\Lambda_h}{\alpha_0\left(\lambda_h^*+\mu\right)-\delta\lambda_h^*};\;\;\;\;\;\;I_h^*=\frac{\lambda^*_h S^*_h}{\alpha_0} 
=\frac{\lambda_h^*\Lambda_h}{\alpha_0\left(\lambda_h^*+\mu\right)-\delta\lambda_h^*},
$$
and
$$
S_v^*=\frac{\Lambda_v}{\eta+\lambda_v^*},\;\;\;\;\,\;I^*_v=\frac{\lambda_v^*\Lambda_v}{\eta(\eta+\lambda^*_v)}
$$
where $\alpha_1=(\alpha+\mu)/(\alpha+\mu+\delta)$, and
$$
\lambda^*_h=\frac{p_1\beta\,I^*_v}{\pi I^*_h + S^*_h},\;\;\;\;\;\;
\lambda^*_v=\frac{\pi p_2\beta\,I^*_h}{\pi I^*_h + S^*_h}.
$$
Therefore,
\begin{equation}\label{lambst}
\lambda_h^*=\frac{p_1\beta\Lambda_v \lambda_v^*\left[\alpha_0\left(\lambda_h^*+\mu\right)-\delta_h\lambda_h^*\right]}{\eta
\left(\pi\lambda_h^*\Lambda_h+\alpha_0\Lambda_h\right)
\left(\eta+\lambda_v^*\right)},
\end{equation}
and
$$
\lambda_v^*=\frac{\pi p_2 \beta \lambda_h^*}{\pi\lambda_h^*+\alpha_0}.
$$
Substituting this last in (\ref{lambst}) one gets the quadratic equation
$$
A_0(\lambda_h^*)^2+B_0\lambda_h^*+C_0 = 0,
$$
where the coefficients are given by (\ref{A0})-(\ref{C0}). Note that $A_0>0$. Then, the thesis follows by applying the Descartes' rule of signs. \#

\medskip

Theorem \ref{exendemic}, point \emph{(ii)}, has established the possibility of multiple equilibria for $R_0<1$. Note that the inequality $B_0^2-4A_0C_0>0$, written in terms of the basic reproduction number, may be written
$$
R_0>R_b,
$$
where
$$
R_b=\frac{\left(\eta\alpha_0\Lambda_h(2\eta+p_2\beta)-p_1p_2\beta^2\Lambda_v(\alpha+\mu)\right)^2-4\eta^3\Lambda_h^2\alpha_o^2(\eta+p_2\beta)}{4\eta^3\Lambda_h^2\alpha_0^2(\eta+p_2\beta)}.
$$
Therefore, model (\ref{model}) has two positive equilibria if 
\begin{equation}\label{condbif}
\max\left\{R_a, R_b\right\} < R_0 < 1.
\end{equation}

\section{Bifurcation analysis}
\label{bifurc}

In this Section we prove that the occurrence of multiple endemic equilibria for $R_0<1$ comes from a backward bifurcation. This will also give information on the local stability of endemic equilibria. To this aim, we study the centre manifold near the criticality (at $E_0$ and $R_0=1$) by using the approach developed in \cite{ccso,duhucc,vawa02}, which is based on the general centre manifold theory \cite{guho}. In short, this approach establishes that the normal form representing the dynamics of the system on the centre manifold is given by 
\begin{equation}
\dot u = a u^2+ b \mu u,
\end{equation}
where,
\begin{equation}
\label{Coeff_a}
a=\frac{{\mathbf v}}{2} \cdot D_{{\mathbf x}{\mathbf x}}{\mathbf f}({\mathbf x}_0,0){\mathbf w}^2 \equiv \frac12 \displaystyle \sum_{k,i,j=1}^{n} v_{k} w_{i} w_{j} \displaystyle \frac{\partial^{2}f_{k}}{\partial x_{i} \partial x_{j}} ({\mathbf x}_0,0),
\end{equation}
and
\begin{equation}
\label{Coeff_b}
b= {\mathbf v} \cdot D_{{\mathbf x}\xi}{\mathbf f}({\mathbf x}_0,0){\mathbf w} \equiv
\displaystyle \sum_{k,i=1}^{n} v_{k} w_{i} \displaystyle \frac{\partial^{2}f_{k}}{\partial x_{i} \partial \xi} ({\mathbf x}_0,0).
\end{equation}

\noindent Note that in the (\ref{Coeff_a}) and (\ref{Coeff_b}) $\xi$ denotes a bifurcation parameter to be chosen, $f_i$'s denote the right hand side of system (\ref{model}), ${\mathbf x}$ denote the state vector, ${\mathbf x}_0$ the disease--free equilibrium $E_0$ and \textbf{v} and \textbf{w} denote the left and right eigenvectors, respectively, corresponding to the null eigenvalue of the Jacobian matrix of system (\ref{model}) evaluated at the criticality.

In our case, let us take $p_2$ as bifurcation parameter. Then, $R_0=1$ is equivalent to
\begin{equation}
\label{p2crit}
p_2=p_2^{crit}:= \frac{\Lambda_h\eta^2\left(\alpha+\eta+\delta\right)}{\pi p_1\mu\Lambda_v\beta^2}.
\end{equation}

It can bee seen that Theorem \ref{stabE0} implies that $b>0$ (see \cite{vawa02}). Therefore, the sign of coefficient (\ref{Coeff_a})  `decides' the direction of the bifurcation occurring at $p_2=p_2^{crit}$.  Precisely, if $a > 0$, then system (\ref{model}) exhibits a backward bifurcation at $R_0=1$. If $a < 0$, then the system exhibits a forward bifurcation at $R_0=1$ \cite{ccso,duhucc,vawa02}.

\begin{theorem}\label{th:bif}
If \begin{equation}
\label{Theta} \Theta :=
\frac{\alpha+\mu}{\Lambda_h}-2\frac{\pi\mu}{\Lambda_h}-\frac{\alpha_0\eta}{\beta p_1 \Lambda_v}>0,
\end{equation}
then system (\ref{model}) exhibits a backward bifurcation at $R_0=1$. If the reversed inequality holds, then the system
exhibits a forward bifurcation at $R_0=1$.
\end{theorem}

\noindent {\bf Proof.} Let us begin by observing that the matrix
\begin{equation}
J(E_0, p_2^{crit})=
\left(\begin{array}{cccc}
-\mu & \delta & 0 & -p_1\beta \\
0 & -\alpha_0 & 0 &  p_1\beta \\
0 & -\frac{\eta\alpha_0}{p_1\beta} & -\eta & 0 \\
0 &  \frac{\eta\alpha_0}{p_1\beta} & 0 & -\eta
\end{array}
\right),
 \end{equation}
admits a simple zero eigenvalue and the other eigenvalues are real and negative. Hence, when $p_2=p_2^{crit}$ (or, equivalently, when $R_{0}=1$),  the disease-free equilibrium $E_{0}$ is a nonhyperbolic equilibrium.

Denote by  ${\bf v}=(v_{1},v_{2},v_{3})$, and ${\bf
w}=(w_{1},w_{2},w_{3})^{T}$, a left and a right eigenvector
associated with the zero eigenvalue, respectively, such that ${\bf
v} \cdot {\bf w}=1$. We get:
$$
{\bf v}=\left(0,\frac{\eta\alpha_0}{p_1\beta(\eta+\alpha_0)},0,\frac{\alpha_0}{\eta+\alpha_0} \right), \;\;\;\;
{\bf w}=\left(-\frac{p_1\beta(\alpha+\mu)}{\mu\alpha_0},\frac{p_1\beta}{\alpha_0},-1,1\right)^{T}.
$$
Taking into account of system (\ref{model}) and considering only the nonzero components of the left eigenvector {\bf v}, it follows that
$$
\begin{array}{ll}
a = & 
      2v_2w_2w_4\displaystyle \frac{\partial^{2}f_{2}}{\partial I_h \partial I_v} (E_0,p_2^{crit})+
      2v_4w_1w_2\displaystyle \frac{\partial^{2}f_{4}}{\partial S_h \partial I_h} (E_0,p_2^{crit})+\\ \\
    & 2v_4w_2w_3\displaystyle \frac{\partial^{2}f_{4}}{\partial I_h \partial S_v} (E_0,p_2^{crit})+
     v_4w_2^2\displaystyle \frac{\partial^{2}f_{4}}{\partial I_h^2} (E_0,p_2^{crit}).
\end{array}
$$
Now it can be checked that
$$
\displaystyle \frac{\partial^{2}f_{2}}{\partial I_h \partial I_v} (E_0,p_2^{crit})=
-\frac{\pi p_1\beta}{S_{h0}},\;\;\;\;\;\;\;\displaystyle \frac{\partial^{2}f_{4}}{\partial S_h \partial I_h} (E_0,p_2^{crit})=
-\frac{\pi p_2^{crit}\beta}{S^2_{h0}}S_{v0},
$$
and
$$
\displaystyle \frac{\partial^{2}f_{4}}{\partial I_h \partial S_v} (E_0,p_2^{crit})=\frac{\pi p_2^{crit}\beta}{S_{h0}},\;\;\;\;\;\;\;
\displaystyle \frac{\partial^{2}f_{4}}{\partial I_h^2} (E_0,p_2^{crit})=
-\frac{2\pi^2p_2^{crit}\beta}{S_{h0}^2}S_{v0},
$$
In view of (\ref{p2crit}) we have
$$
\displaystyle \frac{\partial^{2}f_{4}}{\partial S_h \partial I_h} (E_0,p_2^{crit})=
-\frac{\eta \alpha_0}{\beta p_1 S_{h0}},\;\;\;\;
\displaystyle \frac{\partial^{2}f_{4}}{\partial I_h \partial S_v} (E_0,p_2^{crit})=
\frac{\eta\alpha_0}{\beta p_1S_{v0}},
$$
and
$$
\displaystyle \frac{\partial^{2}f_{4}}{\partial I_h^2} (E_0,p_2^{crit})=
-\frac{2\pi\eta\alpha_0}{p_1\beta S_{h0}},
$$
where $S_{h0}$ and $S_{v0}$ are given in (\ref{E0}). Then, it follows
$$
a= \frac{2p_1\beta\eta}{(\eta+\alpha_0)} \Theta,
$$
where $\Theta$ is defined in (\ref{Theta}). Therefore, system (\ref{model}) exhibits backward or forward bifurcation at
$R_0=1$ according to the sign of $\Theta$. \#

\medskip

\begin{figure}[t]
\centering
\includegraphics[scale=0.45]{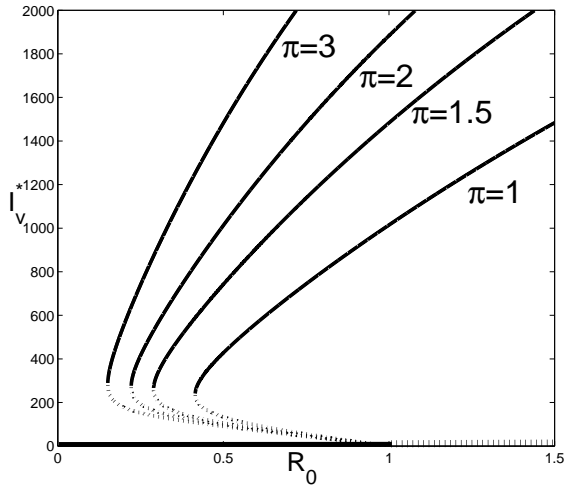}
\caption{\small The backward bifurcation curves in the ($R_0, I_v^*$)--plane as $\pi$ is varied and $b=0.4$. The solid lines represent stability, the dotted lines represent instability. The set of parameter values is given in Table \ref{Tab:param}, with the exception of $p_2$, which has been taken as bifurcation parameter.}
\label{fig:bifdiagP}
\end{figure}

\begin{figure}[t]
\centering
\includegraphics[scale=0.45]{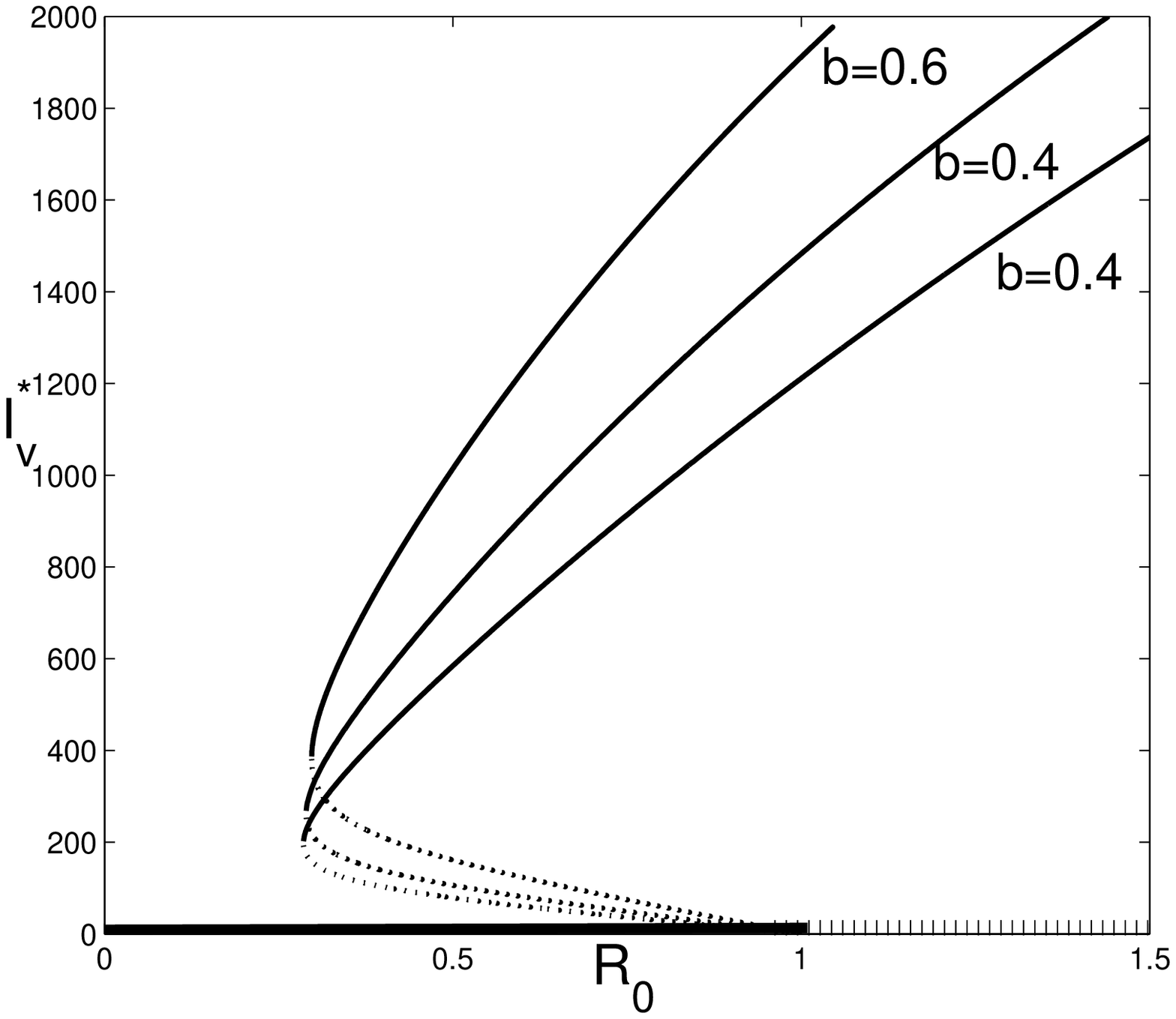}
\caption{\small The backward bifurcation curves in the ($R_0, I_v^*$)--plane as $b$ is varied and $\pi=2$. The solid lines represent stability, the dotted lines represent instability. The set of parameter values is given in Table \ref{Tab:param}, with the exception of $p_2$, which has been taken as bifurcation parameter.}
\label{fig:bifdiagB}
\end{figure}

From (\ref{Theta}) it follows that $\alpha=0$ implies $\Theta<0$, since $\pi\geq 1$. In other words, the disease--induced mortality is responsible for the occurrence of backward bifurcation. However, both the vector--bias parameter $\pi$ and the bed--net usage parameter $b$ have a role in the phenomenon occurrence. This will investigated later on (Section \ref{impact}), where the impact of $\pi$ and $b$ on the quantity $\Theta$ given by (\ref{Theta}) will be investigated. However, this influence may be ``visualized'' in the bifurcation diagram. See Figure \ref{fig:bifdiagP} and Figure \ref{fig:bifdiagB}, where the paths of the endemic states are depicted for various values of $\pi$  and $b$. We use the parameter values in Table \ref{Tab:param}, which are taken from \cite{agetal}. A possible baseline value of the vector--bias parameter is $\pi=2$, since the test performed in \cite{laetal} shows that the number of mosquitoes attracted to gametocytes carrier is the double of those attracted to other individuals (uninfected or carrying non transmittable forms of the parasite). Therefore, we estimate that mosquito that bite humans will do it at a probability $2/3$ if the human is infectious (in the sense that it carries gametocytes) and $1/3$ if the human is susceptible, so that $\pi=2$.

From Figure \ref{fig:bifdiagP} it can be seen that $\pi$ has a huge influence on the saddle-node threshold of $R_{0}$, i.e. the value of $R_0$ below which the only stable equilibrium is the disease--free equilibrium. This threshold decreases as $\pi$ increases, so that high values of $\pi$ make the disease eradication more difficult. From Figure \ref{fig:bifdiagB} we see that this threshold remains substantially unchanged by changing $b$. 

\begin{figure}[t]
\centering
\includegraphics[scale=0.43]{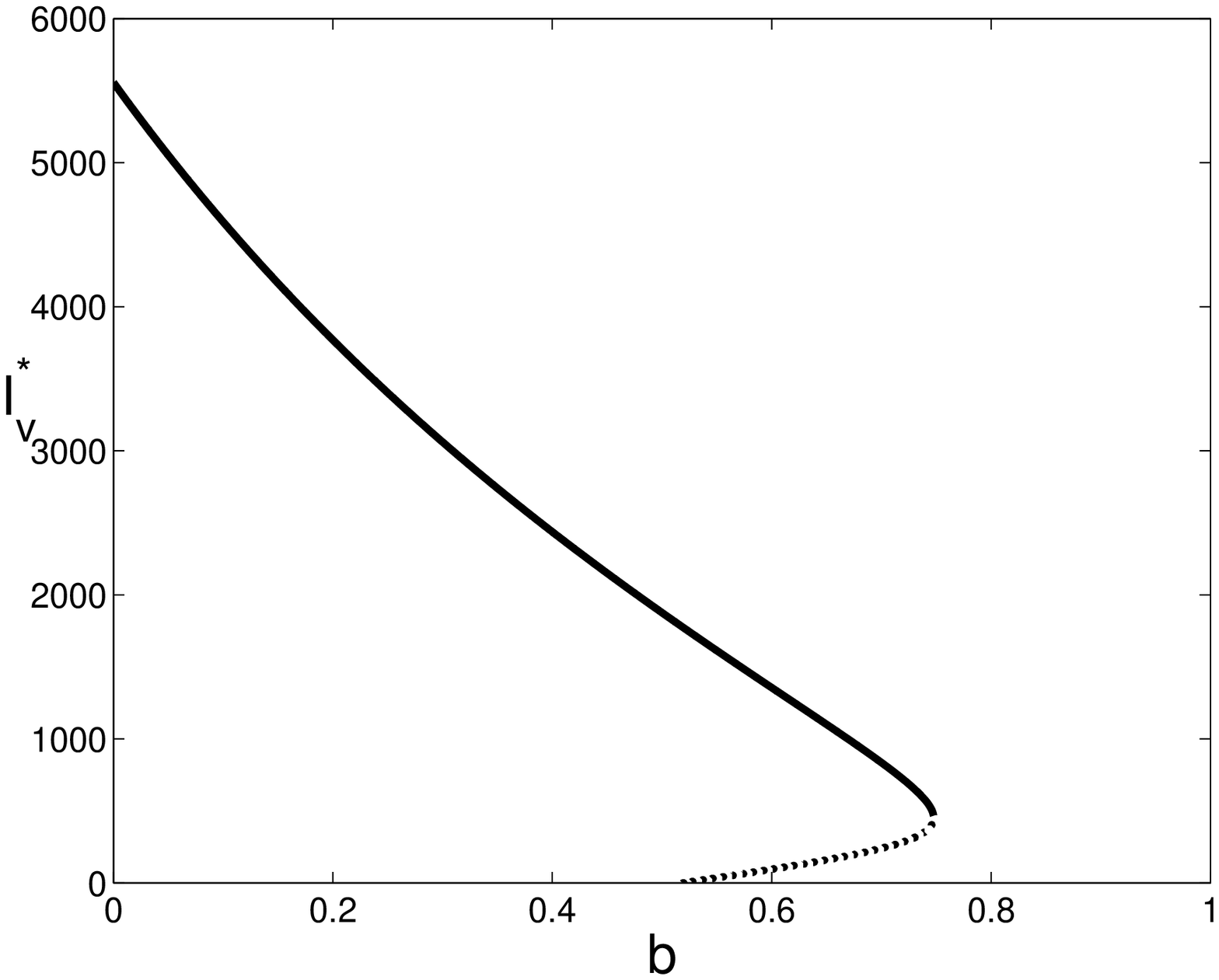} 
\caption{\small The infectious vectors at equilibria as $b$ is varied and $\pi=2$. The solid lines represent stability, the dotted lines represent instability. The set of parameter values is given in Table \ref{Tab:param}, and $p_2=0.6$.}
\label{fig:IvectB}
\end{figure}

Note the bifurcation diagrams are depicted in terms of $R_0$ rather than the bifurcation parameter $p_2$. A consequence of that is the apparent increasing values of the stable branch of $I_v^*$ as $b$ increases (Figure \ref{fig:bifdiagB}). Instead, the values of the stable infectious vectors at equilibrium are decreasing with $b$, as shown in Figure \ref{fig:IvectB}.

\section{Global stability of the endemic equilibrium}
\label{global}

In \cite{agetal} the global stability analysis for the malaria model incorporating bed--net usage has been performed for the special case $\alpha=0$ (no disease--induced human deaths), which implies that the total human population converges to $S_{h0}$. If we do the same here, model (\ref{model}) may be simplified by assuming that the human population is at equilibrium, i. e. $N_h(t)=S_{h0}$, for all $t>0$. In this case, the second generation approach given in \cite{ccetal,dihe} may be fruitfully used to establish that the disease--free equilibrium is globally asymptotically stable when $R_0<1$. This has been done in \cite{agetal} and could be analogously obtained here for model (\ref{model}). Instead, here we focus on the endemic states. 

In the previous sections we have established that if $R_0>1$, then there exists an unique endemic equilibrium, say $E$, for system
(\ref{model}). We now  prove that such an equilibrium is globally asymptotically stable in the interior of the feasible region $\Omega$. This means that the disease becomes endemic for $R_0>1$ regardless of the number of initial cases in both the human and vector populations. Furthermore, this result preclude the possibility  that $E$ destabilizes via onset of oscillations, as it may happen when the human behavior is included in epidemic modeling. For example, when the human behavior is influenced by the available information on the present and the past disease prevalence \cite{bbaml,bbjmaa,domasa}.

We will use  the geometric approach to global stability due to M. Li and J. Muldowney \cite{limu}, which is briefly summarized in
the appendix. The essential of the method is that several sufficient conditions are required for the  global stability of
$E$. Precisely:\\
\emph{(i)} the uniqueness of $E$ in the interior of the set $\Omega$ (i.e. condition (H.1) in the appendix);\\
\emph{(ii)} the existence of an absorbing compact set in the interior of $\Omega$ (i.e. condition (H.2));\\
\emph{(iii)} the fulfillment of a Bendixson criterion (i.e. inequality (\ref{q2})).

Proving the fulfillment of these conditions for a four--dimensional system, like (\ref{model}), can be done but the procedure becomes particularly involved (see for example \cite{baetal,bbctJBD,guetal}). However, the assumption of a constant total population (vector or humans) allows to reduce model (\ref{model}) to a more tractable three-dimensional system. In order to avoid the restriction $\alpha=0$, we assume that the total population at equilibrium is that of vectors, instead of humans, as done in \cite{agetal}, i. e. we assume $N_v(t)=V$ (const.), for all $t>0$. From (\ref{TotPop}) and (\ref{E0}) it immediately follows that $V=S_{v0}$. Under this assumption, from (\ref{model}) we get
\begin{equation}\label{Ridotto}
\begin{array}{ll}
\dot S_{h}=&\Lambda_h-\lambda_h(b,\pi)S_h-\mu S_{h}+\delta I_{h}\\
\dot I_{h}=&\lambda_h(b,\pi)S_h-\left(\alpha+\mu+\delta\right) I_{h}\\
\dot I_{v}=& \lambda_v(b,\pi)\left(V-I_v\right)-\eta(b) I_{v},
\end{array}
\end{equation}
where $\lambda_h(b,\pi)$ and $\lambda_v(b,\pi)$ are given in (\ref{lambda}). This system may be studied in the feasible region
$$
\Omega_0=\left\{(S_h,I_h,I_v)\in{\bf R}^3: 0\leq N_h(t) \leq\frac{\Lambda_h}{\mu},\;\;0\leq I_v(t)\leq V\right\}.
$$
We have the following result: 

\begin{theorem} \label{GASGM} Suppose that in system (\ref{model}) the vector population is at equilibrium. If $R_0>1$, then the unique endemic equilibrium $E$ of (\ref{model}) is globally asymptotically stable.
\end{theorem}

\noindent{\bf Proof.} When $R_0>1$, system (\ref{Ridotto}) satisfies conditions (H.1)-(H.2). In fact, the existence and uniqueness of $E$ has been shown in Section \ref{sec:model}. On the other hand, the instability of $E_0$ (Theorem \ref{stabE0}),  implies the uniform persistence  \cite{frruta}, i.e. there exists a constant $c>0$ such that any solution $(S_h(t), I_h(t), I_v(t))$ with $(S_h(0), I_h(0), I_v(0))$ in the interior of $\Omega_0$, satisfies
  $$\min \left\{ \liminf_{t\rightarrow\infty}S_h(t), \liminf_{t\rightarrow\infty}I_h(t), \liminf_{t\rightarrow\infty}I_v(t)\right\}>c.$$
The uniform persistence together with boundedness of $\Omega_0$, is
equivalent to the existence of a compact set $K$ in the interior of
$\Omega_0$ which is absorbing for (\ref{model}), see \cite{husc}.
Thus, (H.1) is verified. Moreover, $E$ is the only equilibrium in
the interior of $\Omega_0$, so that (H.2) is  also verified. \\ It
remains to find conditions for which the Bendixson criterion given
by (\ref{q2}) is verified.

To this aim, note first that the Jacobian matrix corresponding to system (\ref{Ridotto}) is given by
$$
J=\left(\begin{array}{ccc}
    -\frac{\partial \lambda_h}{\partial S_h}S_h-\lambda_h-\mu & -\frac{\partial \lambda_h}{\partial I_h}S_h+\delta & -\frac{\partial \lambda_h}{\partial I_v}S_h \\ \\
    \frac{\partial \lambda_h}{\partial S_h}S_h+\lambda_h & \frac{\partial \lambda_h}{\partial I_h}S_h-(\alpha+\mu+\delta) & \frac{\partial \lambda_h}{\partial I_v}S_h \\ \\
    \frac{\partial \lambda_v}{\partial S_h}\left(V-I_v\right) & \frac{\partial \lambda_v}{\partial I_h}\left(V-I_v\right) &
    -\lambda_v-\eta
\end{array}\right).
$$
From this we get the second additive compound matrix
$$
J^{[2]}\left( S_{h},I_{h},I_{v}\right)=
\begin{pmatrix}
-a_{11} & \frac{\partial \lambda_h}{\partial I_v}S_h & \frac{\partial \lambda_h}{\partial I_v}S_h\ \\ \ \\
\frac{\partial \lambda_v}{\partial I_h}\left(V-I_v\right)  & -a_{22} & -\frac{\partial \lambda_h}{\partial I_h}S_h+\delta\ \\ \ \\
 -\frac{\partial \lambda_v}{\partial S_h}\left(V-I_v\right)  & \frac{\partial \lambda_h}{\partial S_h}S_h+\lambda_h& -a_{33}
 \end{pmatrix},
$$
where,
$$
\begin{array}{lllllllll}
a_{11}&= \alpha+2\mu+\delta+\lambda_h+\left(\frac{\partial \lambda_h}{\partial S_h}-\frac{\partial \lambda_h}{\partial I_h}\right)S_h,\\ \\
a_{22}&= \mu+\eta+\lambda_h+\lambda_v+\frac{\partial \lambda_h}{\partial S_h}S_h,\\ \\
a_{33}&= \alpha+\mu+\delta+\eta+\lambda_v-\frac{\partial \lambda_h}{\partial I_h} S_h.
\end{array}
$$
Choose now the matrix $P=P(S_{h},I_{h},I_{v})=diag(1,I_{h}/I_{v},I_{h}/I_{v})$. Then $P_{f}P^{-1}=diag(0,\; \dot I_{h}/I_{h}-\dot I_{v}/I_{v},\; \dot I_{h}/I_{h}-\dot I_{v}/I_{v})$, and the matrix $B=P_{f}P^{-1}+PJ^{[2]}P^{-1}$ can be written in block form as
\[
B =\left[
\begin{array}{cc}
B_{11} & B_{12}  \medskip \\
B_{21} & B_{22}
\end{array}
\right],
\]
where,
\begin{align*}
B_{11}&=-\alpha-2\mu-\delta-\lambda_h-\left(\frac{\partial \lambda_h}{\partial S_h}-\frac{\partial \lambda_h}{\partial I_h}\right)S_h,
\end{align*}
\[
B_{12}=\left[
\begin{array}{cc}
\displaystyle \frac{S_hI_v}{I_h}\frac{\partial \lambda_h}{\partial I_v} &
\displaystyle \frac{S_hI_v}{I_h}\frac{\partial \lambda_h}{\partial I_v}
\end{array}
\right],
\]
\[
B_{21} =\left[
\begin{array}{cc}
\displaystyle\frac{\partial \lambda_v}{\partial I_h}\frac{I_h}{I_v}\left(V-I_v\right)  \medskip \\
-\displaystyle\frac{\partial \lambda_v}{\partial S_h}\frac{I_h}{I_v}\left(V-I_v\right)
\end{array}
\right],
\]
\[
B_{22} =\left[
\begin{array}{cc}
\displaystyle\frac{\dot I_{h}}{I_{h}}-\frac{\dot I_v}{I_{v}}-a_{22}&  -\displaystyle\frac{\partial \lambda_h}{\partial I_h}S_h+\delta \\  \\
\displaystyle\frac{\partial \lambda_h}{\partial S_h}S_h+\lambda_h & \displaystyle\frac{\dot I_{h}}{I_{h}}-\frac{\dot I_v}{I_{v}}-a_{33}
\end{array}
\right].
\]
Choose now the vector norm $|\cdot|$ in $\bf{R}_{+}^{3}$ given by
\begin{align*}
|(x,y,z)|&=\max\{|x|,|y|+|z|\}.
\end{align*}
Let $\sigma(\cdot)$ denote the Lozinski\u{i} measure with respect
to this norm. Using the method of estimating $\sigma(\cdot)$  in
\cite{limu}, we have
$$\sigma(B)\leq \sup\left\{ g_1, g_2 \right\} := \sup\left\{ \sigma_1(B_{11})+|B_{12}|, ~\sigma_1(B_{22})+|B_{21}| \right\},
$$
where $|B_{21}|$, $|B_{12}|$ are matrix norms  with respect to the
$L^1$ vector norm and $\sigma_1$ denotes the Lozinski\u{\i}
measure with respect to the $L^1$ norm\footnote{i.e., for the
generic matrix $A=(a_{ij})$, $|A|=\max_{1\leq k\leq
n}\sum_{j=1}^n|a_{jk}|$ and $\mu (A)=\max_{1\leq k\leq n}(a_{kk}+
{\sum_{ {j=1(j\neq k)}}^n|a_{jk}|})$.}. Since $B_{11}$ is scalar,
its Lozinski\u{i} measure with respect to any norm in
$\bf{R}_{+}$  is equal to $B_{11}$. Therefore,
\begin{equation*}
\sigma_{1}(B_{11})=-\alpha-2\mu-\delta-\lambda_h-\left(\frac{\partial \lambda_h}{\partial S_h}-\frac{\partial \lambda_h}{\partial I_h}\right)S_h.
\end{equation*}
Moreover,
\begin{align*}
\sigma_{1}(B_{22})&=\max\left\{\frac{\dot I_{h}}{I_{h}}-\frac{\dot I_{v}}{I_{v}}-a_{22}+\frac{\partial \lambda_h}{\partial S_h}S_h+\lambda_h,\; \frac{\dot I_{h}}{I_{h}}-\frac{\dot I_v}{I_{v}}-a_{33}-\frac{\partial \lambda_h}{\partial I_h}S_h+\delta\right\}\\
&=\max\left\{\frac{\dot I_{h}}{I_{h}}-\frac{\dot I_{v}}{I_{v}}-\mu-\eta-\lambda_v,\; \frac{\dot I_{h}}{I_{h}}-\frac{\dot I_v}{I_{v}}-a_{33}-\alpha-\mu-\eta-\lambda_v\right\}\\
&=\frac{\dot I_{h}}{I_{h}}-\frac{\dot I_{v}}{I_{v}}-\mu-\eta-\lambda_v,
\end{align*}
and
\begin{equation*}
|B_{12}|=\displaystyle \frac{S_hI_v}{I_h}\frac{\partial \lambda_h}{\partial I_v},\;\;\;\;\;
|B_{21}|=\displaystyle\left(\frac{\partial \lambda_v}{\partial I_h}-\frac{\partial \lambda_v}{\partial S_h}\right)\frac{I_h}{I_v}\left(V-I_v\right).
\end{equation*}
Therefore,
\begin{equation}\label{g1}
g_{1}=-\alpha-2\mu-\delta-\lambda_h-\left(\frac{\partial \lambda_h}{\partial S_h}-\frac{\partial \lambda_h}{\partial I_h}\right)S_h+
\displaystyle \frac{S_hI_v}{I_h}\frac{\partial \lambda_h}{\partial I_v},
\end{equation}
and
\begin{equation}\label{g2}
g_{2}=\displaystyle \frac{\dot I_{h}}{I_{h}}-\frac{\dot I_{v}}{I_{v}}-\mu-\eta-\lambda_v+\displaystyle\left(\frac{\partial \lambda_v}{\partial I_h}-\frac{\partial \lambda_v}{\partial S_h}\right)\frac{I_h}{I_v}\left(V-I_v\right).
\end{equation}
From (\ref{Ridotto}) we get
\begin{equation}\label{1g1}
\lambda_h\frac{S_{h}}{I_{h}}=\frac{\dot I_h}{I_{h}}+(\alpha+\mu+\delta),
\end{equation}
and
\begin{equation}\label{1g2}
\frac{\lambda_v}{I_v}\left(V-I_{v}\right)=\frac{\dot I_v}{I_{v}}+\eta.
\end{equation}
Observe that:
\begin{equation}
\displaystyle \frac{S_hI_v}{I_h}\frac{\partial \lambda_h}{\partial I_v}=\lambda_h\frac{S_{h}}{I_{h}},
\end{equation}
and substitute (\ref{1g1}) into (\ref{g1}) and (\ref{1g2}) into (\ref{g2}), to get
$$
g_{1}=\displaystyle \frac{\dot I_{h}}{I_{h}}-\mu-\lambda_h-\left(\frac{\partial \lambda_h}{\partial S_h}-\frac{\partial \lambda_h}{\partial I_h}\right)S_h,
$$
and
$$
g_{2}=\displaystyle \frac{\dot I_{h}}{I_{h}}-\mu-\lambda_v+\left[\displaystyle\left(\frac{\partial \lambda_v}{\partial I_h}-\frac{\partial \lambda_v}{\partial S_h}\right)\frac{I_h}{I_v}-\frac{\lambda_v}{I_v}\right]\left(V-I_v\right).
$$
Now, taking into account of (\ref{lambda}), we have
\begin{equation}\label{g11}
g_{1}=\displaystyle \frac{\dot I_{h}}{I_{h}}-\mu-\frac{p_1\beta\pi I_vI_h}{(\pi I_h+S_h)^2},
\end{equation}
and
\begin{equation}\label{g21}
g_2=\displaystyle \frac{\dot I_{h}}{I_{h}}-\mu-\frac{\pi p_2\beta I_h}{\pi I_h+S_h}-
\frac{\pi(\pi-1)p_2\beta I_h}{\left(\pi I_h + S_h\right)^2}\left(V-I_v\right).
\end{equation}
Equalities (\ref{g11}) and (\ref{g21})  imply
\begin{align*}
\sigma(B)\leq \frac{\dot I_h}{I_{h}}-\mu\label{Sup}.
\end{align*}
Along each solution $(S_{h}(t),I_{h}(t),I_{v}(t))$ to
(\ref{Ridotto}) with $(S_{h}(0),I_{h}(0),I_{v}(0))\in K$, where $K$ is the compact absorbing set, we have, for $t>T_{0}$,
\begin{align*}
    \frac{1}{t}\int^{t}_{0}\sigma(B)ds\leq\frac{1}{t}\int^{T_{0}}_{0}\sigma(B)ds+\frac{1}{t}\ln\frac{I_{h}(t)}{I_{h}(T_{0})}-\mu\frac{t-T_{0}}{t},
\end{align*}
which implies $\overline{q}_{2}<-\mu/2 <0$, where ${q}_{2}$ is given by (\ref{q2}), so that the proof is completed. \#

\section{Role played by parameters $b$ and $\pi$ in model dynamics}
\label{impact}

We begin by discussing the impact of $b$ and $\pi$ on the basic reproduction number $R_0$. Observe that from (\ref{R0}), taking into account of (\ref{etab}), it can be checked that
$$
\frac{dR_0}{db}=\frac{2\pi p_1p_2\mu\Lambda_v\left[\beta(b)\beta'(b)\eta^2(b)-\eta(b)\eta'(b)\beta^2(b)\right]}{\Lambda_h(\alpha+\mu+\delta)\eta^4(b)}<0
$$
and
$$
\frac{dR_0}{d\pi}=\frac{p_1p_2\mu\Lambda_v\beta^2(b)}{\Lambda_h\eta^2(b)(\alpha+\mu+\delta)}>0
$$
This means, as expected, that bed--net usage is beneficial, in the sense that an increase of bed--net usage produces a reduction of the basic reproduction number, whereas increasing the attractiveness of infected humans to mosquitoes produces an increasing of $R_0$. As shown in the previous sections, the minimal value of $R_0$ below which the infection cannot maintain itself in the population (at least for small perturbation of the disease--free equilibrium) depends on if or not the bifurcation at $R_0=1$ is subcritical (backward) or supercritical (forward).

In case of forward bifurcation, this minimal value if the classical threshold $R_0=1$. Therefore, from condition $R_0<1$ we can find the minimal value $b_{crit}$ of bed--net usage ensuring the potential eradication of the disease. From $R_0<1$, in view of (\ref{R0}) and (\ref{etab}), we get
$$
b > \frac{\sqrt{\pi}b_{\max}-\eta_{nat}\sqrt{\tilde\varphi}}{\sqrt{\pi}\left(b_{\max}-b_{\min}\right)+\eta_{bn}\sqrt{\tilde\varphi}} := b_{crit},
$$
where
$$
\tilde\varphi=\frac{(\alpha+\mu+\delta)\Lambda_h}{p_1p_2\mu\Lambda_v}.
$$
When $\pi=1$, we obtain the critical value in absence of vector--bias, say
$$
b_{1}:=\frac{b_{\max}-\eta_{nat}\sqrt{\tilde\varphi}}{\left(b_{\max}-b_{\min}\right)+\eta_{bn}\sqrt{\tilde\varphi}},
$$
which was found in \cite{agetal}. Being $b_{crit}>b_{1}$, for $\pi>1$, it can be deduced that this critical value of bed--net usage increases when the mosquitoes preference for biting infected humans is taken into account. 

An analogous approach may be employed to assess the role of $b$ and $\pi$ on the minimal value of $R_0$ necessary to avoid endemic states in case of backward bifurcation (see condition (\ref{condbif}))

Now, we want to assess which of the two parameters $b$ and $\pi$ has the greatest influence on changes of $R_0$ values and hence the greatest effect in determining whether the disease may be cleared in the population. To this aim, we provide a local sensitivity analysis of the basic reproduction number (see e.g. \cite{chetal}).\\
Denote by $\Psi$ the generic parameter of system \eqref{model}. We evaluate the \emph{normalised sensitivity index} 
$$
S_{\Psi}=\frac{\Psi}{R_0}\frac{\partial R_0}{\partial \Psi},
$$
which indicates how sensitive $R_0$ is to changes of parameter $\Psi$. A positive (resp. negative) index indicates that an increase in the parameter value results in an increase (resp. decrease) in the $R_0$ value. In our case, we have:
\[ 
S_{\pi}=\frac{\pi}{R_0}\frac{\partial R_0}{\partial \pi}=1,
\]
and
\[ 
S_{b}=\frac{b}{R_0}\frac{\partial R_0}{\partial b}=
-2b \left(\frac{\left(\beta_{\max}-\beta_{\min}\right)}{\beta}+\frac{\eta_{bn}}{\eta}\right).
\]
The quantity $S_b$ is negative, as expected, but its magnitude depends on the parameter values. For example, choosing the values in Table \ref{Tab:param} and $\pi=2$, it can be checked that it is decreasing with $b$ and $S_b\approx -1$ when $b\approx 0.24$. This means that when the mosquito attraction to gametocytes carriers is the double of attraction to uninfected individuals, $R_0$ is most sensitive to changes in the bed--net usage only when the rate of adopters is up to 24\%. 

\begin{figure}[t]
\begin{center}
\begin{tabular}{cc}
$
\begin{array}{c}
\mbox{\includegraphics[scale=0.48]{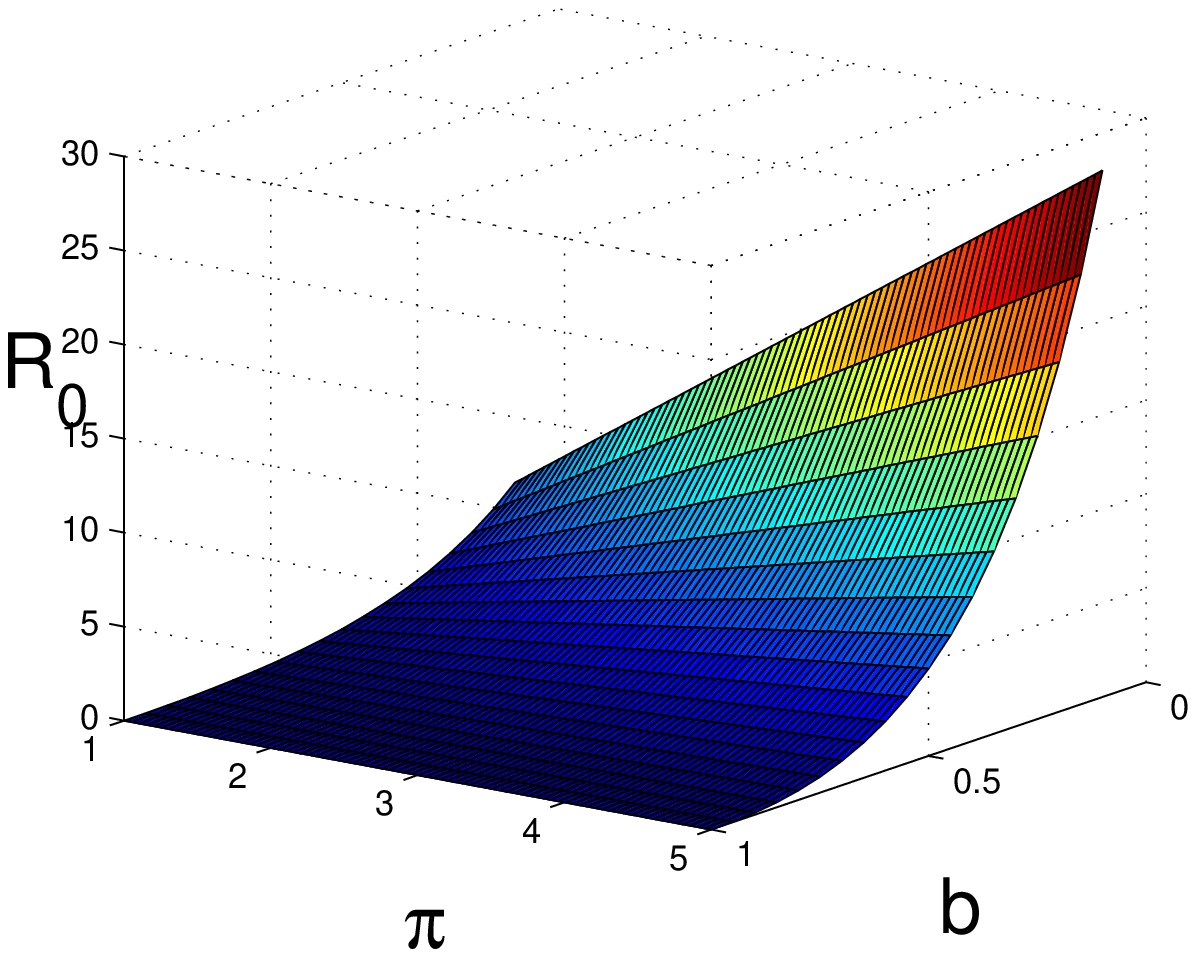}} \\ 
\end{array}
$
$
\begin{array}{c}
\mbox{\includegraphics[scale=0.48]{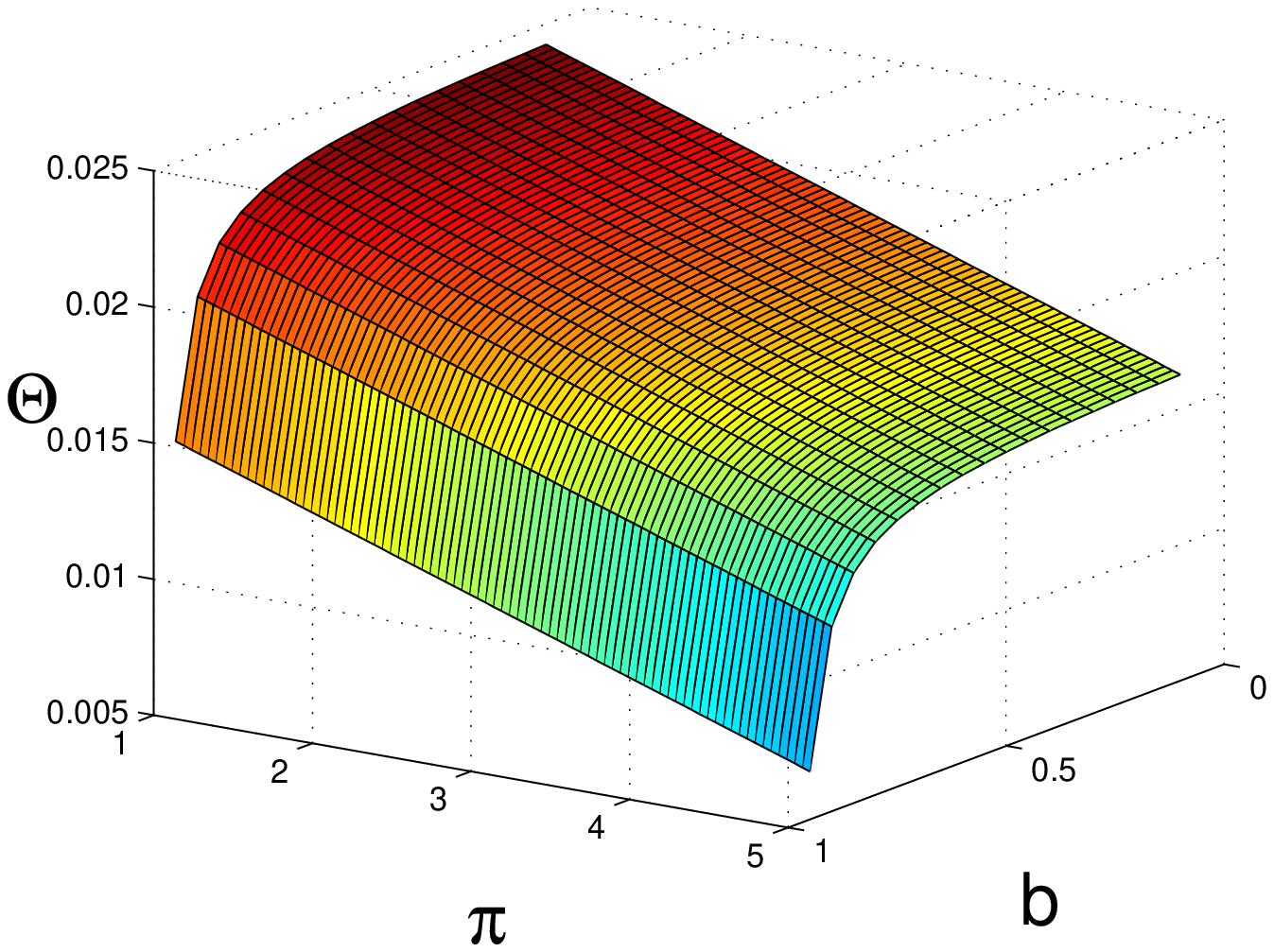}}\\ 
\end{array}
$
\end{tabular}
\caption{\small The basic reproduction number $R_0$ (\emph{left}) and the quantity $\Theta$, given by (\ref{Theta}), plotted as functions of the vector--bias parameter $\pi$ and the bed net usage $b$. The set of parameter values is given in Table \ref{Tab:param}.}
\label{fig:theta}
\end{center}
\end{figure}

The variation of $R_0$ by changing $b$ and $\pi$ can be seen in Figure \ref{fig:theta} (left panel). It is evident the harmful result of a high mosquito preference to infectious together with a low bed-net usage.

We conclude this section by showing the effect of  both the vector--bias and the bed--net usage on the backward bifurcation occurrence. This can be done by analyzing the response of parameter $\Theta$ of changing $\pi$  and $b$. It is easy to check that
$$
\frac{d \Theta}{d\pi}=-\frac{2\mu}{\Lambda_h}<0,
$$
and
$$\frac{d \Theta}{db}=-\frac{\alpha_0}{p_1\Lambda_v}\left[\frac{\eta_{bn}}{\beta_{\max}-b\left(\beta_{\max}-\beta_{\min}\right)}+
\frac{\left(\eta_{nat}+b \eta_{nb}\right)\left(\beta_{\max}-\beta_{\min}\right)}{\left(\beta_{\max}-b\left(\beta_{\max}-\beta_{\min}\right)\right)^2}
\right]<0.
$$
Therefore, even if the system undergoes a backward bifurcation, the phenomenon gets very little observable when the bed--net usage and vector--bias parameters are high enough. This two parameters affects $\Theta$ in different ways, as shown in Figure  \ref{fig:theta} (right panel). Clearly, $\Theta$  reduces linearly with $\pi$, whereas low--middle values of $b$ does not impact too much $\Theta$, which decreases much more rapidly when $b$ is near to its maximum $b=1$. This means that small variation of very high values of bed--net usage produces a big reduction of $\Theta$, so that the phenomenon becomes, in fact, negligible.

\section{Conclusions}
\label{concl}

As far as we know, the simultaneous effects on malaria transmission of \emph{vector--bias} (i. e. the enhanced attractiveness of infectious humans to mosquitoes) and human behavior (here represented by bed--net usage), has never been studied before. We propose a theoretical approach based on modeling and analysis of Mathematical Epidemiology. 

The baseline model used here is the vector--bias malaria model considered in \cite{bbcruz2}, which has been extended by adopting the bed--net usage modeling proposed in \cite{agetal}. The ``merging'' of this two modeling approaches allows to assess that  mosquitoes preference for biting infected humans increases the minimum level of bed--net usage necessary to potentially eradicate  the disease. Moreover, it is shown that low bed-net usage or high vector--bias produces an increase of $R_0$ and therefore the shift towards the globally stable endemicity. We have also seen that the backward bifurcation, and hence the existence of multiple endemic states under the classical threshold $R_0=1$, is a phenomenon essentially due to the disease--induced death rate of humans. Therefore, this phenomenon can be observed only in regions where this rate is very high. However, both $b$ and $\pi$ contributes to reduce its relevance, since their increase reduces the value of $\Theta$ in Theorem \ref{th:bif}, although only variation of very high values of bed--net usage impacts $\Theta$ in a relevant way.

Furthermore, compared to \cite{agetal} and \cite{bbcruz2}, we have the following aspect of novelties:
\begin{itemize}
\item With respect to \cite{agetal}, we have proved that the occurrence of multiple endemic equilibria for $R_0<1$ comes from a backward bifurcation. Using the bifurcation analysis, we are also able to get information on the local stability of the endemic equilibrium emerging from the bifurcation. Furthermore, we have proved that the endemic equilibrium is globally asymptotically stable for $R_0>1$.  This result is new, since in \cite{agetal} the global stability analysis has been performed only for the disease--free equilibrium, and for the special case $\alpha=0$. We needed the assumption of a total vector population. A similar assumption was done in \cite{agetal} to analyse the stability of the disease--free equilibrium.
\item Compared with \cite{bbcruz2}, we consider immigration of both humans and vectors. This means that the first integral given by the constant total vector population does not longer hold. As a consequence, the bifurcation analysis is performed for a four dimensional model (instead of a three dimensional one). Furthermore, the global stability of the endemic equilibrium is proved in terms of the generic force of infection. Hence, the validity of the result may be easily checked for any form of the force of infection. 
\end{itemize}

In conclusion, we stress that the insecticide--treated bed--nets (ITNs) are a non--pharmaceutical intervention to control malaria, which in principle  is used by humans independently of their infectious status. On the other hand, the enhanced attractiveness of infectious humans to mosquitoes is a phenomenon related to host manipulation by malaria parasite and therefore is specifically related to the infectious status. The interplay of this two aspects and the impact on malaria transmission is not immediate at a glance. We found that encouraging bed--net usage is an effective way to control malaria because it reduces the contact rate and, in turn, this reduces the basic reproduction number and may avoid the occurrence of sub--threshold endemic states. However, our analysis shows that mosquitoes preference for biting infected humans may negatively impact the response of malaria dynamics to bed-net usage. These considerations are the result of a theoretical approach. Real data, when available, could validate our findings.


\appendix

\section{The geometric method to global stability}

We deem appropriate to recall the geometric approach to global stability of steady states as developed by Li
and Muldowney \cite{limu}. Consider the autonomous dynamical
system
\begin{equation}
\label{ode} \dot{x}=f(x),
\end{equation}
where $f:D\rightarrow {\bf R}^n$, $D\subset {\bf R}^n$ open set and simply connected and $f\in C^1(D)$. Let $x^*$ be an equilibrium of (\ref{ode}), i.e. $f(x^*)=0$. We recall that $x^*$ is said to be {\it globally stable} in $D$ if it is locally stable and all trajectories in $D$ converge to $x^*$.\\
The following theorem holds \cite{limu}:
\begin{theorem}
\label{thlimu} Assume that:\\
($H1$) there exists a compact absorbing set $K\subset D$;\\
($H2$) the equation (\ref{ode}) has a unique equilibrium $x^*$ in $D$. Then
$x^*$ is globally asymptotically stable in $D$ provided that a
function P(x) and a Lozinski\u{\i} measure ${\cal L}$ exist such
that the inequality
\begin{equation}
\label{q2} 
\overline{q}_2:=\limsup_{t\rightarrow \infty} \sup_{x_0\in \Omega}
{\displaystyle \frac{1}{t} \int_0^t {\cal L}(B(x(s,x_0)))ds}<0,
\end{equation}
is satisfied.
\end{theorem}
\noindent In (\ref{q2}) the quantity $B$ is given by
\[
B=P_fP^{-1}+PJ^{[2]}P^{-1},
\]
where $P(x)$ be a {\small $( \begin{array}{c} n \\ 2
\end{array} ) \times ( \begin{array}{c} n \\ 2 \end{array} )$}
matrix-valued function that is $C^1$ on $D$, and the matrix $P_f$ is
\[
(p_{ij}(x))_f=(\partial p_{ij}(x)/\partial x)^T \cdot f(x)=\nabla
p_{ij}\cdot f(x).
\]
Furthermore, $J^{[2]}$ is the second additive compound matrix of
the Jacobian matrix $J$, i.e. $J(x)=Df(x)$. Finally,  ${\cal L}$ indicates the Lozinski\u{\i} measure of $B$ with respect
to a vector norm $\left|\cdot\right|$ in ${\bf R}^N,~N=(
\begin{array}{c} n \\ 2 \end{array} )$ (see \cite{ma})
\[
{\cal L}(B)={\displaystyle \lim_{h\rightarrow 0^+} \frac{\vline~ I+hB
~\vline-1}{h}}.
\]
We note that for a $n \times n$ Jacobian matrix $J=(J_{ij})$, $J^{[2]}$ is a {\small $(
\begin{array}{c} n \\ 2 \end{array}) \times ( \begin{array}{c} n
\\ 2 \end{array} )$} matrix (for a survey on compound matrices and
their relations to differential equations see \cite{mu}) and in
the special case $n=3$, one has
\[
J^{[2]}= \left[
\begin{array}{ccc}
J_{11}+J_{22} &  J_{23} & -J_{13} \\
 J_{32} & J_{11}+J_{33} & J_{12} \\
 -J_{31} & J_{21} & J_{22}+J_{33}
\end{array}
\right].
\]


\end{linenumbers}
\end{document}